\documentclass[11pt]{article}
\usepackage{amssymb}

\def\sgn{\hbox{sgn}}

\def\intR{\int_{\R}}

\newtheorem{theorem}{Theorem}[section]
\newtheorem{pro}{Proposition}[section]
\newtheorem{lem}{Lemma}[section]
\newtheorem{cor}{Corollary}[section]

\newtheorem{rem}{Remark}[section]
\newcommand{\R}{{ I\!\!R}}

\newcommand{\ep}{\varepsilon}

\def\supp{\mathop{\rm supp\,}\nolimits}
\def\qed{$\hfill \square$}

\newcommand{\re}[1]{(\ref{#1})}

\begin{document}
\begin{center}
\noindent {\Large \bf{Stability of multi antipeakon-peakons profile}}
\end{center}
\vskip0.2cm
\begin{center}
\noindent
{\bf Khaled El Dika$^\sharp$ and  Luc Molinet$^\clubsuit$}\\
{\small
$\sharp$ L.A.G.A., Institut Galil\'ee, Universit\'e Paris-Nord,\\
93430 Villetaneuse, France.\\
$\clubsuit$  L.M.P.T., UFR Sciences et Techniques, Universit\'e de Tours, Parc Grandmont, 37200 Tours, FRANCE.} \vskip0.3cm
 \noindent
% E-mail :
khaled@math.univ-paris13.fr\\
% E-mail :
Luc.Molinet@lmpt.univ-tours.fr
\end{center}
\vskip0.5cm \noindent {\bf Abstract.} { The Camassa-Holm equation
possesses well-known peaked solitary waves that can travel to
both directions. The positive ones travel to the right and are
called peakon whereas the negative ones travel to the left and are
 called antipeakons.
 Their orbital stability  has been established by
 Constantin and Strauss in \cite{CS1}. In \cite{EL2} we have proven the stability of trains of peakons.
  Here, we continue  this study by extending the stability result to the case of
 ordered  trains of anti-peakons and peakons.}\vspace*{4mm} \\

%%%%%%%% INTRO
\section{Introduction}

\noindent The  Camassa-Holm equation (C-H),
\begin{equation}
u_t -u_{txx}=- 3 u u_x +2 u_x u_{xx} + u u_{xxx}, \quad
(t,x)\in\R^2, \label{CH}
\end{equation}
can be derived as a model for the propagation of unidirectional
shalow water waves over a flat bottom  by  writing the
Green-Naghdi equations in Lie-Poisson Hamiltonian form and then
making an asymptotic expansion which keeps the Hamiltonian
structure (\cite{CH1}, \cite{Johnson}). Note that the Green-Naghdi equations arise as approximations to the  governing equations
 for  shallow-water medium-amplitude regime which captures more nonlinear effects than the classical 
  shallow-water small amplitude KdV regime and thus can accommodate models for breaking waves
   (cf. \cite{AL}, \cite{CL}, \cite{CE}). The Camassa-Holm equation was also found independently by Dai
\cite{dai} as a model for nonlinear waves in cylindrical
hyperelastic rods and was, actually, first discovered by the
method of recursive operator by Fokas and Fuchsteiner \cite{FF}
as an example of bi-Hamiltonian equation. Let us also mention that it has also a geometric derivation as a 
 re-expression of geodesic flow on the diffeomorphism group on the line (cf. \cite{K1}, \cite{K2}) and that this 
  framework is instrumental in showing that the Least Action Principle holds for this equation (cf. \cite{C},
   \cite{CK}).

(C-H) is completely integrable (see \cite{CH1},\cite{CH2},
\cite{C1} and \cite{CGI}). It possesses among  others the
following invariants
\begin{equation}
 E(v)=\int_{\R} v^2(x)+v^2_x(x)
\, dx \mbox{ and } F(v)=\int_{\R} v^3(x)+v(x)v^2_x(x) \, dx\;
\label{E}
\end{equation}
and can be written in Hamiltonian form as
\begin{equation}
\partial_t E'(u) =-\partial_x F'(u) \quad .
\end{equation}
Camassa and Holm \cite{CH1} exhibited peaked solitary waves
solutions to (C-H) that  are
 given by
$$ u(t,x)=\varphi_c(x-ct)=c\varphi(x-ct)=ce^{-|x-ct|},\; c\in\R.$$
They are called peakon whenever $ c>0 $ and antipeakon  whenever
$c<0$.  Let us point out here   that the feature of the peakons that their profile is smooth, except at the crest 
 where it is continuous but the lateral tangents differ, is similar to that of the waves of greatest height, i.e. traveling waves of largest possible amplitude which are solutions to the governing equations for water waves (cf. 
 \cite{C2}, \cite{CE2} and \cite{T}).  
Note that (C-H) has to be rewriten as
\begin{equation}
u_t +u u_x +(1-\partial_x^2)^{-1}\partial_x (u^2+u_x^2/2)=0
\label{CH3} \; .
\end{equation}
to give a meaning to these solutions. Their stability seems not
to enter the general framework developed
 for instance in \cite{Benjamin}, \cite{GSS}. However, Constantin and
Strauss \cite{CS1} succeeded in proving their orbital stability by
a direct approach. In \cite{EL2} we combined  the general strategy
initiated in \cite{MMT}(note that due to  the reasons mentioned
above, the general method of \cite{MMT} is not directly
applicable here ), a monotonicity result proved in \cite{EL} on
the part of the energy $ E(\cdot) $ at the right of a localized
solution traveling to the right and  localized versions of the
estimates established in \cite{CS1} to derive the stability of
ordered trains of peakons. In this work we pursue this
 study by proving the stability of ordered trains of anti-peakons and peakons. The main new ingredient is a
  monotonicity result on the part of the functional $ E(\cdot)-\lambda F(\cdot) $, $ \lambda\ge 0 $,  at the right
   of a localized solution traveling
   to the right. It is worth noticing that  the sign of $ \lambda $ plays a crucial role in our
   analysis.

Before stating the main result let us  introduce the function
space where we will
 define the flow of the equation. For $ I $  a finite or infinite interval of $ \R $,
  we denote by $ Y(I)  $ the function space\footnote{$ W^{1,1}(\R) $ is the space of $ L^1(\R) $
  functions with derivatives in $ L^1(\R) $ and $ BV(\R) $ is the space of function with bounded variation}
  \begin{equation}
Y(I):= \Bigl\{ u\in C(I;H^1(\R)) \cap L^\infty(I;W^{1,1}(\R)), \;
u_x\in L^\infty(I; BV(\R))\Bigr\} \, . \label{theoweak}
\end{equation}
In \cite{CM1}, \cite{dan1} and \cite{E} (see also \cite{L}) the
following existence and uniqueness result for this class of initial data  is
 derived.
\begin{theorem} \label{wellposedness}
Let $ u_0\in H^1(\R) $ with $ m_0:=u_0-u_{0,xx} \in {\cal M}(\R)
$ then there exists $ T=T(\|m_0\|_{\cal M})>0 $ and a unique
solution $ u \in Y([-T,T]) $ to (C-H) with initial data $ u_0 $.
The functionals $ E(\cdot) $ and $ F(\cdot) $ are constant along
the trajectory and if $ m_0 $ is such that  \footnote{ ${\mathcal
M}(\R)$ is the space of Radon measures on $ \R $ with bounded
total variation. For $ m_0\in {\mathcal M}(\R)$ we denote
respectively by $ m_0^- $ and $ m_0^+ $ its positive and negative
part.}
 $ \supp m_0^- \subset ]-\infty, x_0] $
and $ \supp m_0^+ \subset [x_0,+\infty[ $ for some $ x_0\in \R $ then $ u $  exists for all positive times and belongs to $ Y([0,T]) $ for all $ T>0 $. \\
Moreover,  let $ \{u_{0,n}\}\subset H^1(\R) $  such that $  u_{0,n}\to u_0 $ in $ H^1(\R) $  with $
\{m_{0,n}:=u_{0,n}-\partial^2_x u_{0,n} \} $  bounded in $ {\cal M} (\R) $, 
$ \supp m_{0,n}^- \subset ]-\infty, x_{0,n}] $
and 
$ \supp m_{0,n}^+ \subset [x_{0,n},+\infty[ $ for some sequence $ \{x_{0,n}\}\subset \R  $. Then, for all $ T> 0
$,
\begin{equation}
u_n \longrightarrow u \mbox{ in } C([0,T]; H^1(\R)) \; .
\end{equation}
\end{theorem}
 Let
 us emphasize that the global existence result when the negative part of $ m_0 $ lies completely to the left of its positive
   part is proven in \cite{E} and  that the last assertion of the above theorem is not explicitly contained in this paper.
    However, following the same arguments as those developed in these works (see for instance Section 5 of \cite{L}), one can prove that there exists a subsequence
  $ \{u_{n_k}\} $ of solutions of \re{CH} that converges in $ C([0,T]; H^1(\R)) $ to some solution $ v$ of \re{CH} belonging to $ Y([0,T[) $. Since $ u_{0,n_k} $ converges to
   $ u_0 $ in $ H^1 $, it follows that $ v(0)=u_0 $ and thus $ v=u $ by uniqueness. This
    ensures that the whole sequence $ \{u_n\} $ converges to $ u $ in $C([0,T]; H^1(\R))$
     and concludes the proof of the last assertion.
   \begin{rem}
   It is worth pointing out that recently, in \cite{BC1} and \cite{BC2}, Bressan and Constantin
 have constructed global conservative and dissipative solutions of the Camassa-Holm
  equation for any initial data in $ H^1(\R) $. However, even if for the conservative solutions, $ E(\cdot) $ and $ F(\cdot) $
   are conserved quantities, these  solutions are not known to be continuous with values in $ H^1(\R) $. Therefore even one single peakon
   is not known to be orbitally stable
  in this class of solutions. For this reason we will work in the class
    of solutions constructed in Theorem \ref{wellposedness}.
   \end{rem}

We are now ready to state our main result.
\begin{theorem} \label{mult-peaks}
Let be given $ N $ non vanishing velocities
$c_1<..<c_{k}<0<c_{k+1}<..< c_N $ . There exist  $ \gamma_0 $, $
A>0 $, $ L_0>0 $
 and $ \varepsilon_0>0 $ such that if $ u \in Y([0,T[)$, with $ 0<T\le \infty $,
   is a solution of (C-H) with initial data $ u_0 $ satisfying
 \begin{equation}
 \|u_0-\sum_{j=1}^N \varphi_{c_j}(\cdot-z_j^0) \|_{H^1} \le \varepsilon^2 \label{ini}
 \end{equation}
 for some  $ 0<\varepsilon<\varepsilon_0$ and $ z_j^0-z_{j-1}^0\ge L$,
with $ L>L_0 $, then there exist $x_1(t), ..,x_N(t) $ such that
\begin{equation}
\sup_{[0,T[} \|u(t,\cdot)-\sum_{j=1}^N
\varphi_{c_j}(\cdot-x_j(t)) \|_{H^1} \le
A(\sqrt{\varepsilon}+L^{-{1/8}})\;  .\label{ini2}
\end{equation}
Moreover there exists $ C^1 $-functions  $ {\tilde x_1}, ..,
{\tilde x_N} $ such that, $ \forall j\in \{1,..,N\} $,
\begin{equation}\label{gy}
|x_j(t)-{\tilde x_{j}}(t)|=O(1) \mbox{ and }
 \frac{d}{d t} {\tilde x_j}(t) = c_j +O(\varepsilon^{1/4}) +O(L^{-\frac{1}{16}}),  \forall t\in[0,T[ \; .
\end{equation}
 \end{theorem}
\begin{rem} \label{rem1}
We do not know how to prove the monotonicity result in Lemma
\ref{monotonicitylem}, and thus Theorem \ref{mult-peaks}, for
solutions that are only in $ C([0,T[;H^1(\R)) $ which is the
hypothesis required for the stability of a single peakon (cf.
\cite{CS1}). Note anyway that there exists no well-posedness
result in the class  $ C([0,T[;H^1(\R)) $ for general initial
data in $ H^1(\R) $.  On the other hand, according to Theorem
\ref{wellposedness} above, $ u\in Y([0,T[) $ as soon as $u_0 \in
H^1(\R) $ and $ (1-\partial_x^2) u_0 $ is a Radon measure with
bounded variations. 
\end{rem}
\begin{rem}
Note that under the hypotheses of Theorem \ref{mult-peaks}, $$\sum_{j=1}^N \varphi_{c_j}(\cdot-z_j^0)$$
 belongs to the class  $ v\in H^1(\R) $with $m:=v-v_{xx} \in {\cal M}(\R) $, 
 $ \supp m^- \subset ]-\infty, x_0] $
and $ \supp m^+ \subset [x_0,+\infty[ $ for some $ x_0\in \R $. 
Therefore, in view of Theorem
\ref{wellposedness}, Theorem \ref{mult-peaks} leads to the orbital stability (for positive times) of such ordered sum of antipeakons and 
 peakons with respect to $ H^1 $-perturbations that keep the initial data in this same class.
\end{rem}
As discovered by Camassa and Holm \cite{CH1}, (C-H) possesses also special
solutions called multipeakons given by
$$
u(t,x)=\sum_{i=1}^N p_j(t) e^{-|x-q_j(t)|} \,,
$$
where $ (p_j(t),q_j(t)) $ satisfy the differential system :
\begin{equation} \label{systH}
\Bigl\{
\begin{array}{l}
\dot{q}_i=\sum_{j=1}^N p_j e^{-|q_i-q_j|} \\
\dot{p}_i = \sum_{j=1}^N p_i p_j \sgn (q_i-q_j) e^{-|q_i-q_j|}\; .
\end{array}
\Bigr.
\end{equation}
 In \cite{Beals} (see also \cite{Beals0} and \cite{CH1}), the  asymptotic behavior of the   multipeakons is
  studied. In particular, the limits as $ t $ tends to $ +\infty $ and $ -\infty $
   of $ p_i(t) $ and $ \dot{q_i}(t) $ are determined.
 Combining these asymptotics  with the  preceding theorem  and the continuity with respect to 
  initial data stated in Theorem \ref{wellposedness}  we  get the following result
on the stability for positive times of the variety $ {\cal N}_{N,k} $ of $ H^1(\R) $
defined for $ N\ge 1 $ and $ 0\le k \le N $  by
$$
{\cal N}_{N,k}:= \Bigl\{ v=\sum_{i=1}^N p_j e^{-|\cdot-q_j|}, \,
(p_1,..,p_{N})\in (\R_-^*)^k\times (\R_+^*)^{N-k} ,\,
q_1<q_2<..<q_N \,\Bigr\} \; .
$$
\begin{cor} \label{cor-mult-peaks}
Let be given $ k $  negative real numbers $ p_1^0,.., p_{k}^0 $, $
N-k $ positive real numbers $ p_{k+1}^0,..,p_N^0 $ and $ N $ real
numbers $ q_1^0< ..< q_N^0 $. For any $ B> 0 $ and any $ \gamma
>0 $ there exists $ \alpha>0 $ such that if $ u_0\in H^1(\R) $ is
such that $ m_0:=u_0-u_{0,xx} \in
{\mathcal M}(\R) $ with $ \supp m_0^- \subset ]-\infty, x_0] $
and $ \supp m_0^+ \subset [x_0,+\infty[ $ for some $ x_0\in \R $,
and satisfies
\begin{equation}
\|m_0\|_{\cal M}\le B \quad  \mbox{ and }\quad
\|u_0-\sum_{j=1}^N p_j^0 \exp (\cdot-q_j^0) \|_{H^1}\le \alpha
\label{ini3}
\end{equation}
 then
\begin{equation}
\forall t\in\R_+, \quad \inf_{P\in \R_-^k\times\R_+^{N-k},Q\in \R^N}
\|u(t,\cdot)-\sum_{j=1}^N p_j \exp (\cdot-q_j) \|_{H^1} \le
\gamma\; . \label{ini33}
\end{equation}
Moreover, there exists $ T>0 $ such that
\begin{equation}
\forall t\ge T, \quad \inf_{Q\in {\mathcal G}} \|u(t,\cdot)-\sum_{j=1}^N \lambda_j \,\exp
(\cdot-q_j) \|_{H^1} \le \gamma \label{ini4}
\end{equation}
where $ {\mathcal G}:=\{Q\in \R^N, \, q_1<q_2<..<q_N\} $ and $
\lambda_1<..<\lambda_N $ are the eigenvalues of the  matrix
 $ \Bigl( p_j^0 e^{-|q_i^0-q_j^0|/2}\Bigr)_{1\le i,j\le N}
$.
\end{cor}
\begin{rem}
Again, note that for $(p_1^0,..,p_N^0)\in (\R_-^*)^k\times(\R_+^*)^{N-k}$ and 
$  q_1^0<..<q_N^0 $,  $$\sum_{j=1}^N p_j^0 \exp (\cdot-q_j^0)$$ belongs to the class 
 $ v\in H^1(\R) $with $m:=v-v_{xx} \in {\cal M}(\R) $, 
 $ \supp m^- \subset ]-\infty, x_0] $
and $ \supp m^+ \subset [x_0,+\infty[ $ for some $ x_0\in \R $. Corollary \ref{cor-mult-peaks} thus ensures that 
 the variety $ {\cal N}_{N,k} $ is stable with respect to $ H^1$-perturbations that keeps the initial data
  in this same class.
\end{rem}
This paper is organized as follows. In the next section we sketch the main points of 
the proof of Theorem ~\ref{mult-peaks} whereas the complete proof is given in  Section~\ref{sec-multi}. After having  controlled  the distance
between the different bumps of the solution we establish
 the new monotonicity result and state  local
versions of  estimates involved in the stability of a single
peakon.
 Finally, the proof of  Theorem~\ref{mult-peaks} is completed in Subsection \ref{end-proof}.
%%%%%%%%%%%%%%%%%%%
%%%%%%%%%%%% One peakon

\section{Sketch of the proof}\label{section1peak}
Our proof as in \cite{MMT} combined the stability of a single peakon and a monotonicity result for   functionals related to the conservation laws.
Recall that the stability proof  of  Constantin and
Strauss (cf. \cite{CS1}) is principally based on
the following lemma of $\cite{CS1}$.
\begin{lem}\label{1peakon-lemme}
For any $u\in H^1(\R)$, $ c\in \R $  and $\xi\in\R$,
\begin{equation}\label{eq1}
E(u)-E(\varphi_c)= \|u-\varphi_c(\cdot-\xi)\|^2_{H^1}+4c(u(\xi)-c).
\end{equation}
For any $u\in H^1(\R)$, let $M=\max_{x\in\R}\{u(x)\}$, then
\begin{equation}\label{eq2}
M E(u)-F(u)\geq \frac{2}{3}M^3.
\end{equation}
\end{lem}
Indeed, with this lemma at hand, let $u\in C([0, T[ ; H^1(\R))$ be a solution
of $\re{CH}$ with
$\|u(0)-\varphi_c\|_{H^1}\leqslant\varepsilon^2$ and let
$\xi(t)\in\R$ be such that $u(t,\xi(t))=\max_{\R} u(t,\cdot) $.
 Assuming that $ u(t)$ is sufficiently  $ H^1 $-close to $ \{r\in \R, \varphi(\cdot-r)\} $,  setting
 $\delta=c-u(t,\xi(t))$, and using that $ E(u(t))=E(u_0)=2c^2+O(\varepsilon^2) $ and
 $ F(u(t))=F(u_0)=\frac{4}{3}c^3+O(\varepsilon^2) $, \re{eq2} leads to 
 $$
  \delta^2 (c-\delta/3)\le O(\varepsilon^2))\Longrightarrow \delta \lesssim \varepsilon
 $$
  and then \re{eq1} yields
 \begin{equation}\label{eq33}
 \|u(t)-\varphi_c(\cdot-\xi(t))\|_{H^1}\lesssim \sqrt{\varepsilon} \,.
 \end{equation}
 This proves the stability result. At this point, a crucial remark  is that, instead of using
 the conservation of $ E $ and $F $, we can only use that, for any fixed
  $ \lambda  \ge 0 $, $ E(\cdot)-\lambda F(\cdot) $ is non increasing. Indeed, for
   $M=\max_{x\in\R}\{u(t,x)\}=u(t,\xi(t))$ and $ \lambda=1/M $,  \re{eq2} then  implies
   $$
   M E(u_0)-F(u_0)\ge \frac{2}{3} M^3
   $$
   and, for $ \lambda=0 $, \re{eq1}  implies
   $$
   E(u_0)-E(\varphi_c)\ge  \|u-\varphi_c(\cdot-\xi(t))\|^2_{H^1}+4c(u(\xi(t))-c).
   $$
   This leads to \re{eq33}  exactly as above. \\
   Now, in \cite{EL2} it is established that \re{eq1} and \re{eq2} almost still hold if one replaces $ E(\cdot) $ and
      $ F(\cdot) $ by their localized version, $ E_j(\cdot) $ and $ F_j(\cdot) $, around the  j$th $ bump. Therefore to prove our result it will somehow suffices
       to prove  that the functionals $  E_j(\cdot) +\lambda  F_j(\cdot) $ are almost decreasing.

  One of the very important discovering of  the works of Martel-Merle is that for  one dimensional dispersive equations with a linear group that
    travels to the left,  the part of the energy at the right of a
      localized solution traveling to the right is almost decreasing.
       In \cite{MMT}
     it is noticed that this holds also for the
     part of the energy at the right of each bump for solutions that are
     close to the sum of solitary waves traveling to the right.
       In this paper we will use that,  for a fixed $ \lambda\ge 0$ and $ j\ge k+1 $,
        if we call by $ I_j =\sum_{q=j}^N ( E_q-\lambda F_q)$ the part of the functionals
      $ E(\cdot)-\lambda F(\cdot) $ that is at the right of the  (j-1)$th
      $ bumps,  then $ I_j(\cdot) $ is almost decreasing in time.
    Since $ I_N=E_N-\lambda F_N $,  we infer from above that the $ Nth $ bump of the solution stays $H^1 $-close to a translation of
       $\varphi_{c_N} $.  Then, since $ I_{N-1} =E_{N-1}-\lambda F_{N-1}+I_N$ and $I_{N-1} $ is almost decreasing,  we obtain that
        $ E_{N-1}-\lambda F_{N-1} $ is also almost decreasing which leads to
        the stability result for the $ (N-1)th $ bump. Iterating
         this process until $ j=k+1 $,
         we obtain that each bump moving to the right remains
         close to the orbit of the suitable peakon.  Finally, since (C-H) is invariant
          by the change of unknown  $
u(t,x)\to -u(t,-x) $, this also ensures that  each bump moving to
the left remains
         close to  the orbit of the suitable antipeakon. This leads to the
         desired result since the total energy is conserved.

         Actually we will not proceed exactly that way since by
         using such iterative process one loses some power of
         $\varepsilon$ at each step. More precisely this iterative
         scheme would prove Theorem \ref{mult-peaks} but with $
         \varepsilon^\beta $ with $ \beta =4^{1/2-\max(q, N-q)} $
          instead of $\varepsilon^{1/2} $ in \re{ini2}. To derive
         the
         desired  power of $\varepsilon$  we
         will rather sum  all the contributions of bumps that are traveling in the
          same direction and use Abel's summation
         argument to get the stability of all these bumps in the same time.

%%%%%%%%%%%%%%%%%%%%%%%%

%%%%%%%%%%%% Proof of the result

\section{Stability of multipeakons}\label{sec-multi}
%Proof of Theorem~$\ref{mult-peaks}$}

For $ \alpha>0 $ and $ L>0 $ we define the following neighborhood of all the sums
  of N antipeakons and peakons of speed $ c_1,..,c_N $ with spatial shifts $ x_j $ that satisfied
$ x_j-x_{j-1}\ge L $.
  \begin{equation}
U(\alpha,L) = \Bigl\{
u\in H^1(\R), \, \inf_{x_j-x_{j-1}> L} \|u-\sum_{j=1}^N \varphi_{c_j} (\cdot-x_j) \|_{H^1} < \alpha \Bigr\}\; .
\end{equation}

 By the continuity of the map $ t\mapsto u(t) $ from $ [0,T[ $ into $ H^1(\R) $,
to prove the first part of Theorem \ref{mult-peaks} it suffices to prove that there exist
$ A>0 $, $ \varepsilon_0>0 $
and $ L_0>0 $ such  that $ \forall L>L_0 $ and $ 0<\varepsilon<\varepsilon_0 $, if
$u_0$ satisfies \re{ini} and if  for some $ 0<t_0< T $,
\begin{equation}\label{e11}
 u(t)\in U\left(A(\sqrt{\varepsilon}+L^{-{1/8}}),L/2\right)  \textrm{ for all }t\in[0,t_0]
\end{equation}
then
 \begin{equation}\label{e12}
 u(t_0) \in U\left(\frac{A}{2}(\sqrt{\varepsilon}+L^{-{1/8}}),\frac{2L}{3}\right).
\end{equation}
Therefore, in the sequel of this section we will assume \re{e11} for some
$ 0<\varepsilon<\varepsilon_0 $ and $ L>L_0 $, with $ A$, $ \varepsilon_0 $ and $ L_0 $ to be specified later, and we will prove
 \re{e12}.
%%%%%%%%%%%%%%%%%%%%%%%
%%%%%%%%%%%%% Modulation
\subsection{Control of the distance between the peakons}
In this subsection we want to prove that the different bumps of $ u $  that are individualy close to a peakon or an antipeakon  get away  from each others as time is increasing.
 This is crucial in our analysis since we do not know how to manage strong interactions.
  The following lemma is  principally proven in \cite{EL2}.
\begin{lem}\label{eloignement}
Let $ u_0 $ satisfying \re{ini}. There exist $\alpha_0>0$,
$L_0>0$ and $C_0>0$  such that for all  $0<\alpha<\alpha_0$ and
$0<L_0<L$ if $u\in U(\alpha, L/2) $ on $ [0,t_0] $ for some $
0<t_0< T $ then there exist $ C^1
$-functions  $ {\tilde x_1}, .., {\tilde x_N} $ defined on $
[0,t_0] $ such that $\forall t\in  [0,t_0] $,
  \begin{equation}
  \frac{d}{d t} {\tilde x_i}(t) = c_i +O(\sqrt{\alpha}) +O(L^{-1}) , \; i=1,..,N \, ,
  \label{vitesse}
  \end{equation}
  \begin{equation} \label{distH1}
  \|u(t)-\sum_{i=1}^N \varphi_{c_i} (\cdot -{\tilde x_i}(t)) \|_{H^1} =
  O(\sqrt{\alpha}) \, ,
  \end{equation}
\begin{equation}
{\tilde x_i}(t)-\tilde{x}_{i-1}(t) \ge 3L/4+(c_{i}-c_{i-1}) t/2 ,
 \quad i=2,..,N . \label{eloi}
\end{equation}
Moreover, for $  i=1,..,N $, 
   it holds
\begin{equation}
 |x_i(t)-{\tilde x_i}(t)| =O(1)  , \label{prox}
  \end{equation}
  where $ x_{i}(t)\in [\tilde{x}_i(t)-L/4, \tilde{x}_{i}(t)+L/4] $ is any point such that
\begin{equation}
 |u(t,x_i(t))|=\max_{[\tilde{x}_i(t)-L/4, \tilde{x}_{i}(t)+L/4]} |u(t)| .  \label{maxi}
  \end{equation}
\end{lem}
{\it Proof. } We only sketch the proof and refer to \cite{EL2} for details. The strategy is to use a modulation argument to  construct $ N $ $C^1$-functions
 $ t\mapsto {\tilde x_i}(t)$, $i=1,..,N $ on $ [0,t_0] $ satisfying the following orthogonality
  conditions :
  \begin{equation}
\int_{\R} \Bigl( u(t,\cdot) -\sum_{j=1}^N \varphi_{c_j}(\cdot-{\tilde x_j}(t)) \Bigr)
 \partial_x \varphi_{c_i} (\cdot -{\tilde x_i}(t)) \, dx = 0 \; . \label{mod2}
\end{equation}
Moreover,  setting
  \begin{equation}\label{ZZ}
R_Z(\cdot)=\sum_{i=1}^N \varphi_{c_i}(\cdot -z_i)
\end{equation}
 for  any $ Z=(z_1,..,z_N)\in \R^N $, one can check that
\begin{equation}
\|u(t)-R_{{\tilde X}(t)}\|_{H^1} \lesssim C_0 \sqrt{\alpha}\, , \quad \forall t\in
[0,t_0] \; . \label{estepsilon}
\end{equation}
To  prove that the speed of $ {\tilde x}_i $ stays close to $ c_i $,
  we set
 $$
 R_j(t)=\varphi_{c_j}(\cdot-{ \tilde x}_j(t)) \mbox{ and }
 v(t)=u(t)-\sum_{i=1}^N R_j(t)=u(t,\cdot)-R_{{\tilde X}(t)} \; .
$$
and differentiate \re{mod2} with respect to time to get
$$
\int_{\R} v_t  \partial_x R_i =\dot{\tilde x}_i \, \langle
 \partial_x^2 R_i \, ,\,  v  \rangle_{H^{-1}, H^1} \, ,
$$
and thus
\begin{equation}
\Bigl|\int_{\R} v_t  \partial_x R_i\Bigr|
\le |\dot{\tilde x}_i| O(\|v\|_{H^1}) \le |\dot{\tilde x}_i-c_i|
 O(\|v\|_{H^1})+O(\|v\|_{H^1})\; . \label{huhu}
\end{equation}
Substituting $ u $ by $ v+\sum_{j=1}^N R_j $ in \re{CH3}  and
using that $ R_j $
 satisfies
 $$
 \partial_t R_j +(\dot{\tilde x}_j-c_j)  \partial_x R_j + R_j \partial_x R_j
 +(1-\partial_x^2)^{-1}\partial_x  [R_j^2+(\partial_x R_j)^2/2] = 0 \;,
 $$
 we infer that $ v$ satisfies on $ [0,t_0] $,
  \arraycolsep1pt
 \begin{eqnarray}
 v_t&  - & \sum_{j=1}^N (\dot{\tilde x}_j-c_j)  \partial_x R_j
= -\frac{1}{2}  \partial_x \Bigl[(v+\sum_{j=1}^N
 R_j)^2- \sum_{j=1}^N R_j^2  \Bigr] \nonumber \\
 & &-(1-\partial_x^2)^{-1}\partial_x \Bigl[(v+\sum_{j=1}^N R_j)^2- \sum_{j=1}^N R_j^2
 +\frac{1}{2} (v_x +\sum_{j=1}^N \partial_x R_j)^2 -\frac{1}{2}\sum_{j=1}^N
  (\partial_x R_j)^2\Bigr]\; . \nonumber
 \end{eqnarray}
 \arraycolsep5pt
Taking the $ L^2 $-scalar product with $ \partial_x R_i $, integrating by parts, using the
 decay of $ R_j $ and its first derivative, \re{estepsilon} and \re{huhu}, we find
 \begin{equation}
 |\dot{\tilde x}_i-c_i|\Bigl(\|\partial_x R_i \|_{L^2}^2 +O(\sqrt{\alpha}) \Bigr)
 \le O(\sqrt{\alpha}) +O(e^{ -L/8})\; . \label{fofo}
 \end{equation}
Taking $\alpha_0$ small enough and $ L_0 $ large enough we
  get
$ |\dot{\tilde x}_i-c_i| \le (c_i-c_{i-1})/4 $ and thus,  for all $ 0<\alpha<\alpha_0
$ and $ L\ge L_0>3C_0\varepsilon $, it follows from \re{ini}, \re{estepsilon} and \re{fofo}   that
\begin{equation}
{\tilde x}_j(t)-{\tilde x}_{j-1}(t)>
L-C_0\varepsilon+(c_j-c_{j-1}) t/2 , \quad \forall t\in[0,t_0]
\; . \label{xj-xj-1}
\end{equation}
which yields \re{eloi}.\\
Finally from \re{estepsilon} and the continuous embedding of $
H^1(\R)$ into $ L^\infty(\R) $, we infer that
$$ u(t,x) = R_{{\tilde
X}(t)}(x)+O(\sqrt{\alpha}),  \quad \forall x\in \R \, .
$$
Applying this formula with $ x=\tilde{x}_i$  and taking advantage of
 \re{eloi}, we obtain
$$
|u(t,\tilde{x}_i)|
=|c_i|+O(\sqrt{\alpha})+O(e^{-L/4}) \ge
3|c_i|/4 \; . $$
 On the other hand, for $ x\in [\tilde{x}_i(t)-L/4, \tilde{x}_{i}(t)+L/4]\backslash
]{\tilde x_i}(t)-2,{\tilde x_i}(t)+2[ $, we get
$$ |u(t,x)|\le |c_i|
e^{-2}+O(\sqrt{\alpha})+O(e^{-L/4}) \le |c_i|/2 \; . $$ This
ensures that  $x_i$ belongs to $ [{\tilde x_i}-2,{\tilde
x_i}+2] $.
%%%%%%%%%%%%%%%%%%%%%%%%

%
%%%%%%%%%%%%%%%%%%%%%%%%
%%%%%%%%%%%% Monotonie

\subsection{Monotonicity property}\label{Sectmonotonie}
Thanks to the preceding lemma, for $ \varepsilon_0> 0 $ small enough and $ L_0>0 $ large enough, one can construct  $C^1$-functions $ {\tilde x_1}, .., {\tilde x_N} $ defined on $
[0,t_0] $ such that \re{vitesse}-\re{prox} are satisfied.
In this subsection we state the almost monotonicity of  functionals that are very close
 to the $ E(\cdot) -\lambda F(\cdot) $  at the right  of the $ i $th bump, $ i=k,..,N-1 $ of $ u $.
  The proof follows the same lines as in
  Lemma~4.2 in~\cite{EL} but is more delicate since we have also to deal with the functional
   $ F $. Moreover, $ F $ generates a term ( $J_4 $ in \re{go2}) that we are not able to estimate in a suitable way but which fortunately  is of the good  sign.

 Let $ \Psi $ be a  $ C^\infty $ function such that
 $ 0<\Psi\le 1 $,  $ \Psi'>0 $ on $ \R $, $ |\Psi'''|\le 10 |\Psi'| \mbox{ on } [-1,1]
 $,
 \begin{equation}\label{psipsi}
 \Psi(x)=\left\{ \begin{array}{ll}
 e^{-|x|} & \quad x<-1\\
 1-e^{-|x|}& \quad x>1
 \end{array}
 \right. .
 \end{equation}
 Setting $ \Psi_K=\Psi(\cdot/K) $, we introduce for $ j\in \{q,..,N\}
 $ and $ \lambda\ge 0 $,
 $$
 I_{j,\lambda}(t)=I_{j,\lambda,K}(t,u(t))=  \int_{\R}\Bigl(  (u^2(t)+u_x^2(t)) -\lambda (u^3(t)+u
u_x^2(t))\Bigr) \Psi_{j,K}(t) \, dx\,,
 $$
 where $ \Psi_{j,K}(t,x)=\Psi_K(x-y_j(t)) $ with $ y_j(t)$,
 $ j=k+1,..,N $, defined by
 $$
y_{k+1}(t)=\tilde{x}_{k+1}(0)+c_{k+1} t/2 -L/4
 $$
    \begin{equation}\label{defyi}
\mbox{ and } y_i(t)=\frac{{\tilde x_{i-1}}(t)+{\tilde x_i}(t)}{2},\quad
  i=k+2,..,N.
  \end{equation}
    Finally, we set
   \begin{equation}
\sigma_0=\frac{1}{4} \min \Bigl(c_{k+1},c_{k+2}
-c_{k+1},..,c_N-c_{N-1}\Bigr) \; .
   \end{equation}
\begin{pro}\label{monotonicitylem}
  Let $ u\in Y([0,T[) $ be a solution of (C-H) satisfying \re{distH1} on
   $[0,t_0] $.
 There exist $ \alpha_0>0 $ and $ L_0>0 $ only depending on $ c_{k+1} $ and $c_N$ such that
    if $ 0<\alpha<\alpha_0 $ and $ L\ge L_0 $ then for any $ 4\le K \lesssim L^{1/2} $
     and $ 0\le \lambda\le 2/c_{k+1} $,
    \begin{equation}\label{monotonicityestim}
    I_{j,\lambda,K}(t)-I_{j,\lambda,K}(0)\le O( e^{-\frac{\sigma_0 L }{8K}}) ,
    \quad \forall j\in\{k+1,..,N\}, \; \quad \forall t\in [0,t_0] \; .
    \end{equation}
\end{pro}
{\it Proof. }
Let us assume that $ u $ is smooth since the case $ u\in  Y([0,T[) $ follows by modifying slightly the arguments (see Remark 3.2 of \cite{EM}).
  \begin{lem}
\begin{eqnarray}
\frac{d}{dt}\int_{\R} (u^2+u_x^2) g \, dx&=&\int_{\R}(u^3+4uu_x^2)g^{'} \, dx\nonumber\\
&&\hspace*{-10mm}
-\int_{\R}u^3g^{'''} \, dx- 2 \int_{\R}u h g^{'} \, dx.
\label{go}
\end{eqnarray}
and
\begin{eqnarray}
\frac{d}{dt}\int_{\R} (u^3+uu_x^2) g \, dx&=&\int_{\R}(u^4/4+u^2 u_x^2)g^{'} \, dx\nonumber\\
&&\hspace*{-10mm}
+\int_{\R}u^2 h g^{'} \, dx+\int_{\R}( h^2-h_x^2) g^{'} \, dx.
\label{gogo}
\end{eqnarray}
where $ h:=(1-\partial_x^2)^{-1} (u^2+u_x^2/2) $.
\end{lem}
{Proof. }
 Since \re{go} is proven in \cite{EL2}  we concentrate on the proof of \re{gogo}.
\begin{eqnarray}
\frac{d}{dt} \intR (u^3+u u_x^2) g  & = & 3 \intR u_t u^2 g + 2 \intR u_{tx} u_x u g +\intR u_t u_x^2 g
\nonumber\\
& = & 2\intR u_t (u^2+u_x^2/2) g +\intR u_t u^2 g -\intR u_{txx} u^2 g -\intR u_{tx} u^2 g' \nonumber\\
&= &  2\intR u_t (u^2+u_x^2/2) g +\intR (u_t -u_{txx}) u^2 g -\intR u_{tx} u^2 g' \nonumber \\
= I_1+I_2+I_3 \, . \label{hu1}
\end{eqnarray}
Setting $ h:=(1-\partial_x^2)^{-1} (u^2+u_x^2/2) $ and using the equation we get
\begin{eqnarray}
I_1 & = & -2 \intR u u_x (u^2+u_x^2/2) g -2 \intR  g h_x (1-\partial_x^2) h  \nonumber \\
& = & -2 \intR u^3 u_x g -\intR u u_x^3 g -2 \intR h h_x g +2\intR h_x h_{xx} g \nonumber\\
& = & \frac{1}{2} \intR u^4 g' -\intR u u_x^3 g +\intR (h^2-h_x^2) g'  \; .
\end{eqnarray}
In the same way,
\begin{eqnarray}
I_2 & = & -3 \intR u^3 u_x g -\frac{1}{2} \intR \partial_x(u_x^2) u^2 g +\frac{1}{2} \intR \partial^3_x (u^2) u^2 g  \nonumber \\
& = & \frac{3}{4}\intR u^4 g' -\frac{1}{2} \intR \partial_x(u_x^2) u^2 g -\frac{1}{2} \intR
 \partial_x^2 (u^2) \partial_x(u^2) g -\frac{1}{2}\intR \partial_x^2 (u^2) u^2 g'   \nonumber\\
& = & \frac{3}{4}\intR u^4 g'+\intR u u_x^3 g  +\frac{1}{2} \intR  u_x^2 u^2 g' +\frac{1}{4} \intR
 [\partial_x (u^2) ]^2 g' +\intR \partial_x(u^2) u u_x g' \nonumber \\
 & & +\frac{1}{2}\intR \partial_x (u^2) u^2 g^{''}   \nonumber\\
 & = & \frac{3}{4}\intR u^4 g'+\intR u u_x^3 g  +\frac{1}{2} \intR  u_x^2 u^2 g' +\intR
 u^2 u_x^2  g' +2 \intR u^2_x u^2  g' +\intR  u^3 u_x  g^{''}   \nonumber\\
 & = &  \frac{3}{4}\intR u^4 g'-\frac{1}{4}\intR u^4 g^{'''}+\frac{7}{2} \intR  u_x^2 u^2 g'
+\intR u u_x^3 g \; .
\end{eqnarray}
At this stage it is worth noticing that the terms $\intR u u_x^3 g  $ cancels with the one in $ I_1 $.
Finally,
\begin{eqnarray}
I_3 & = & \intR \partial_x ( u u_x) u^2 g' +\intR g' u^2 \partial_x^2 h    \nonumber \\
& = & -2\intR u^2 u_x^2 g' -\intR u^3 u_x g^{''} -\intR u^2( u^2 +u_x^2/2) g' +\intR u^2 h g'    \nonumber\\
& = & -2 \intR u^2 u_x^2 g' +\frac{1}{4} \intR u^4 g^{'''}-\intR u^4 g'  -\frac{1}{2}\intR u^2 u_x^2 g' +\intR u^2 h g'   \nonumber\\
 & = & -\frac{5}{2}  \intR u^2 u_x^2 g' +\frac{1}{4} \intR u^4 g^{'''}-\intR u^4 g' +\intR u^2 h g' \label{hu4}
 \end{eqnarray}
where we used that $  \partial_x^2 (I-\partial_x^2)^{-1} = -I +(I-\partial_x^2)^{-1} $. Gathering \re{hu1}-\re{hu4},  \re{gogo} follows.\qed \hspace*{2mm}\\
Applying \re{go}-\re{gogo} with $ g=\Psi_{j,K} $, $ j\ge k+1 $,  one gets
\begin{eqnarray}
\frac{d}{dt} I_{j,\lambda,K} & : = &
\frac{d}{dt}\int_{\R}  \Psi_{j,K}[(u^2+u_x^2)-\lambda (u^3+ u u_x^2) ] \, dx \nonumber \\
&=&-\dot{y_j}
\int_{\R} \Psi_{j,K}' (u^2+u_x^2)  \nonumber \\
 & & +\int_{\R}\Psi_{j,K}' \Bigl[[(u^3+4uu_x^2)-\lambda\Bigl(\dot{y_j} (u^3+u u_x^2)-(u^4/4 +u^2 u_x^2)\Bigr) \Bigr] \, dx\nonumber\\
&&-\int_{\R}\Psi_{j,K}^{'''} u^3\, dx-\int_{\R}\Psi_{j,K}' (2 u+\lambda u^2) h  \, dx \nonumber \\
&  &- \lambda \intR  \Psi_{j,K}'  (h^2-h_x^2) \, dx\nonumber \\
& = & -\dot{y_j}\int_{\R} \Psi_{j,K}'(u^2+u_x^2) +J_1+J_2+J_3+J_4  \nonumber \\
& \le  & -\frac{c_{k+1}}{2} \int_{\R} \Psi_{j,K}'(u^2+u_x^2) +J_1+J_2+J_3+J_4  \; .
\label{go2}
\end{eqnarray}
We claim that  $ J_4 \le 0 $ and that for  $ i\in \{1,2,3\} $, it holds
\begin{equation}
J_i \le \frac{c_{k+1}}{8} \int_{\R} \Psi_{j,K}' (u^2+u_x^2) + \frac{C}{K} \, e^{-\frac{1}{K}
(\sigma_0 t+L/8)}\; . \label{go4}
\end{equation}
To handle with $ J_1 $ we   divide $ \R $ into two regions $ D_j $ and $ D_j^c $ with
 $$
 D_j=[{\tilde x_{j-1}}(t) +L/4, {\tilde x_j}(t) -L/4]
 $$
 First since from \re{eloi}, for $ x\in D_j^c $ ,
 $$
 |x-y_j(t)| \ge \frac{{\tilde x_{j}}(t)-{\tilde x_{j-1}}(t)}{2}-L/4 \ge \frac{c_j-c_{j-1}}{2} \, t +L/8\, ,
 $$
 we infer from the definition of $ \Psi $ in Section \ref{Sectmonotonie} that
 $$
 \Bigl|\int_{D_j^c} \Psi_{j,K}' \Bigl[[(u^3+4uu_x^2)-\lambda\Bigl(\dot{y_j} (u^3+u u_x^2)-(u^4/4 +u^2 u_x^2)\Bigr) \Bigr] \, dx\Bigr|
 $$
 $$
  \le \frac{C}{K} \, (1+2 \lambda c_N) (\|u_0\|_{H^1}^3+ \|u_0\|_{H^1}^4)
  e^{-\frac{1}{K}(\sigma_0 t +L/8)} \; .
  $$
  On the other hand, on $ D_j $ we notice,  according to \re{distH1}, that
  \begin{eqnarray}
  \|u(t)\|_{L^\infty_{D_j}}& \le &  \sum_{i=1}^N \|\varphi_{c_i}(\cdot -{\tilde x_i(t))\|_{L^\infty}(D_j)}
   + \|u-\sum_{i=1}^N\varphi_{c_i}(\cdot -{\tilde x_i(t)})\|_{L^\infty(D_j)} \nonumber \\
    & \le & C \, e^{-L/8} +O(\sqrt{\alpha}) \; .\label{go3}
    \end{eqnarray}

Therefore, for $  \alpha $ small enough and $ L $ large enough it holds
$$
J_1 \le \frac{c_{k+1}}{8} \int_{\R} \Psi_{j,K}' (u^2+u_x^2) + \frac{C}{K} \, e^{-\frac{1}{K}(\sigma_0 t +L/8)}\; .
$$
Since $ J_2 $  can be handled in exactly the same way, it remains to treat $ J_3 $.
For this, we first notice as above that
  \begin{eqnarray}
& &\hspace*{-15mm} -\int_{D_j^c}(2u+\lambda u^2)  \Psi_{j,K}'
(1-\partial_x^2)^{-1}(u^2+u_x^2/2 ) \nonumber \\
& & \le (2+\lambda\|u\|_{\infty}) \|u\|_{\infty} \sup_{x\in D_j^c}
|\Psi_{j,K}'(x-y_j(t))|\int_{\R} e^{-|x|} \ast (u^2+u_x^2/2 ) \, dx \nonumber \\
 &  & \le \frac{C}{K} \|u_0\|_{H^1}^3 \, e^{-\frac{1}{K}(\sigma_0 t +L/8)}\; ,
 \label{J31}
\end{eqnarray}
since
\begin{equation}
 \forall f\in L^1(\R), \quad (1-\partial_x^2)^{-1} f
  =\frac{1}{2} e^{-|x|} \ast f \; .
 \label{tytu}
 \end{equation}
Now in the region $D_j $, noticing that $ \Psi_{j,K}' $ and
$ u^2+u_x^2/2 $ are non-negative, we  get
 \begin{eqnarray}
 & & \hspace*{-15mm} -\int_{D_j}(2 +\lambda u ) u  \Psi_{j,K}'
(1-\partial_x^2)^{-1}(u^2+u_x^2/2 ) \nonumber \\
 &  \le &
(2+\lambda \|u(t)\|_{L^\infty({D_j})} ) \|u(t)\|_{L^\infty({D_j})} \int_{D_j}\Psi_{j,K}'(
(1-\partial_x^2)^{-1}(2u^2+u_x^2) \nonumber \\
&  \le &  (2+\lambda \|u(t)\|_{L^\infty({D_j})} ) \|u(t)\|_{L^\infty({D_j})} \int_{\R} (2u^2+u_x^2) (1-\partial_x^2)^{-1}
\Psi_{j,K}'\; .
\end{eqnarray}
On the other hand, from the definition of $ \Psi $ in Section \ref{Sectmonotonie}
   and \re{tytu} we  infer  that for $ K\ge 4 $,
  $$
(1-\partial_x^2) \Psi_{j,K}' \ge (1-\frac{10}{K^2}) \Psi_{j,K}' \Rightarrow
(1-\partial_x^2)^{-1} \Psi_{j,K}'\le (1-\frac{10}{K^2})^{-1} \Psi_{j,K}' \; .
  $$
Therefore, taking $ K\ge 4 $ and using \re{go3}   we deduce  for $ \alpha $ small enough and $ L $ large enough that
\begin{equation}
  -\int_{D_j} (2u+\lambda u^2)   \Psi_K'
(1-\partial_x^2)^{-1}(u^2+u_x^2/2)
  \le   \frac{  c_q}{8}
\int_{\R} (u^2+u_x^2/2)
 \Psi_K' \; . \label{J32}
\end{equation}
This completes the proof of \re{go4}. It remains to prove  that $ J_4 $ is non positive.  Recall that
 $  h=(I-\partial_x^2)^{-1} v $ with $ v:=u^2+u_x^2/2 \ge 0 $. Therefore, following \cite{CE1}, it holds
 \begin{eqnarray*}
 h(x)& =& \frac{1}{2} e^{-|\cdot|} \ast v(\cdot) \\
 & = & \frac{1}{2} e^{-x} \int_{-\infty}^x e^y  v(y) \, dy +\frac{1}{2} e^{x} \int_{-\infty}^x e^{-y}  v(y) \, dy
 \end{eqnarray*}
 and
 $$
 h'(x)=-\frac{1}{2} e^{-x} \int_{-\infty}^x e^y  v(y) \, dy +\frac{1}{2} e^{x} \int_{-\infty}^x e^{-y}  v(y) \, dy
 $$
 which clearly ensures that $ h^2 \ge h_x^2 $.  Since $ \Psi_{j,K}'\ge 0$ and $ \lambda \ge 0 $, this
  leads to the  non positivity of $ J_4=- \lambda \intR  \Psi_{j,K}'  (h^2-h_x^2) \, dx $.

 Gathering \re{go2} and \re{go4} we  infer that
$$
\frac{d}{dt}\int_{\R} \Psi_{j,K}[u^2+u_x^2-\lambda(u^3+u u_x^2) ] \, dx\le  -\frac{c_1}{8}
\int_{\R}\Psi_{j,K}' (u^2+u_x^2) +
 \frac{C}{K} (1+\|u_0\|_{H^1}^4) \, e^{-\frac{1}{K}(\sigma_0 t +L/8)} \; .
$$
Integrating this inequality between $ 0 $ and $ t $, \re{monotonicityestim} follows.\\

%%%%%%%%%%%%%%%%%%%
%%%%%%%%%%%% energie

\subsection{Localized estimates}\label{Localized energy estimates}

We define the function $ \Phi_i=\Phi_i(t,x) $, $i=k+1,..,N$,  by $
\Phi_N=\Psi_{N,K}=\Psi_K(\cdot-y_N(t)) $ and for $i=k+1,..,N-1 $
\begin{equation}\label{defphii}
\Phi_i=\Psi_{i,K}-\Psi_{i+1,K}=\Psi_K(\cdot-y_i(t))-\Psi_K(\cdot-y_{i+1}(t))\; ,
\end{equation}
where $ \Psi_{i,K} $ and the $ y_i $'s are defined in Section \ref{Sectmonotonie}. It
is easy to check that the $ \Phi_i$'s are positive functions  and that $\displaystyle \sum_{i=k+1}^N \Phi_{i}\equiv
\Psi_{k+1,K} $.
 We will  take $L/K>4$ so that \re{psipsi} ensures that  $ \Phi_i $ satisfies
  for $i\in \{k+1,..,N\} $, 
 \begin{equation}
 |1-\Phi_{i}| \le  2 e^{-\frac{L}{8K}} \mbox{ on } ] y_i+L/8,  y_{i+1}-L/8[
  \label{de1}
 \end{equation}
 and
\begin{equation}
 |\Phi_{i}| \le 2 e^{-\frac{L}{8K}} \mbox{ on } ]y_i-L/8,
 y_{i+1}+L/8[^c \; ,\label{de2}
 \end{equation}
 where we set $ y_{N+1}:=+\infty $. \\
 It is worth noticing that, somehow, $ \Phi_i(t) $ takes care of only the ith bump of $u(t)$.
 We will use  the following localized  version of $ E $ and $ F $ defined  for
$i\in \{k+1,..,N\}, $ by
 \begin{equation}\label{defEi}
 E_i^t(u) = \int_{\R} \Phi_{i}(t) (u^2+u_x^2) \mbox{ and }
 F_i^t(u)= \int_{\R} \Phi_i(t) (u^3+u u_x^2) \; .
 \end{equation}
{\bf Please note that henceforth we take $K=L^{1/2}/8 $.} \\
The following lemma gives a localized version of \re{eq2}.
Note that the functionals $ E_i $ and $ F_i $ do not depend on
time in the statement below  since we fix $ 
y_{k+1}<..<y_{N+1}=+\infty$.
\begin{lem}\label{m-p-lemme}
 Let  be given $u\in H^1(\R) $ with $\|u\|_{H^1}=\|u_0\|_{H^1} $ and $ N-k $ real numbers $ y_{k+1}<..<y_N$
  with $y_i-y_{i-1} \ge 2L/3 $. For $ i=k+1,..,N $, set $J_i:=]y_i-L/4,y_{i+1}+L/4[ $ with 
  $ y_{N+1}=+\infty$, and assume that  there exist 
  $ x_i\in ]y_i+L/4,y_{i+1}-L/4[ $ such that  $ u(x_i)=\displaystyle \max_{J_i} u:=M_i>0 $.
  Then, defining the functional $ E_i $'s and $ F_i$'s as
  in \re{defphii}-\re{defEi}, 
 it holds
\begin{equation}\label{eq2m}
F_i(u)\leqslant M_i E_i(u)-\frac{2}{3}M_i^3+\|u_0\|_{H^1}^3 O(L^{-{1/2}}), \quad
 i\in \{k+1,..,N\} \, .
\end{equation}
and for any $ x_1<..<x_k$ with $x_k <y_{k+1}-L/4$, setting $X:=(x_{1},..,x_N)\in \R^{N} $, it holds
\begin{equation}\label{eq3m}
 E_i(u)-E(\varphi_{c_i})=E_i (u-R_{X})
 +4c_i (M_i-c_i)+ \|u_0\|_{H^1}^2 O(L^{-1/2})
 , \quad
 i\in \{k+1,..,N\} ,
\end{equation}
where $ R_X $ is defined in \re{ZZ}.
\end{lem}
{\it Proof. } Let $ i\in \{k+1,..,N\}$ be fixed. Following \cite{CS1},
we introduce the function  $g$ defined  by
$$g(x)=\left\{
\begin{array}{l}
u(x)-u_x(x) \; \mbox{ for } \; x<{ x_i} \\
u(x)+u_x(x) \; \mbox{ for }\; x>{ x_i}
\end{array}.
\right.
$$
Integrating by parts we compute
\begin{eqnarray}\label{ug2}
\int ug^2\Phi_{i}&=&\int_{-\infty}^{
x_i}(u^3+uu_x^2-2u^2u_x)\Phi_i
+\int_{ x_i}^{+\infty}(u^3+uu_x^2+2u^2u_x)\Phi_i\nonumber\\
&=&F_i(u)-\frac{4}{3}u({x_i})^3\Phi_i({
x_i})+\frac{2}{3}\int_{-\infty}^{ x_i}u^3\Phi_i^{'}
-\frac{2}{3}\int_{ x_i}^{+\infty}u^3\Phi_i^{'}\, .
\end{eqnarray}
Recall that we take $ K=\sqrt{L}/8 $ and thus  $ |\Phi'|\le C/K = O(L^{-1/2}) $.
 Moreover, since $ x_i\in  ] y_i+L/4,  y_{i+1}-L/4[$, it follows from \re{de1} that     $ \Phi_i({ x
_i})=1+O(e^{-L^{1/2}}) $  and thus
\begin{equation}
\int u g^2 \Phi_i = F_i(u)-\frac{4}{3} M_i^3+\|u\|_{H^1}^3 O(L^{-1/2}) \; .
\end{equation}
On the other hand, with \re{de2} at hand, 
\begin{eqnarray}\label{hg2}
\int u g^2\Phi_i& \le & M_i \int_{J_i} g^2 \Phi_i +  \int_{J_i^c} | u| g^2 \Phi_i \nonumber \\
& \le &  M_i \int_{-\infty}^{+\infty} g^2 \Phi_i + \|u\|_{L^\infty(\R)} \int_{J_i^c}   g^2 \Phi_i \nonumber \\
&\leq &M_i \Bigl( E_i(u)-2\int_{-\infty}^{x_i}uu_x\Phi_i+2\int_{
x_i}^{+\infty}
uu_x\Phi_i \Bigr)+ \|u\|_{H^1}^3 \sup_{x\in J_i^c} |\Phi_i(x)|\nonumber\\
&\leq &M_iE_i(u)-2M_i^3+\|u\|_{H^1}^3 O(L^{-1/2}) \; .
\end{eqnarray}
This proves \re{eq2m}. To prove \re{eq3m}, we
  use the relation between $ \varphi $ an its derivative and
integrate  by parts, to get \arraycolsep2pt
\begin{eqnarray*}
E_i (u-R_X) & = & E_i (u)+E_i (R_X) -2 \int \Phi_i
\Bigl( u\, \varphi_{c_i}
(\cdot-x_i) + u_x \, \partial_x \varphi_{c_i}(\cdot-x_i) \Bigr)  \\
& =& E_i(u)+E_i (R_X) -2 \int \Phi_i u \, \varphi_{c_i}(\cdot-x_i) \\
 & & +2
 \int_{x_i}^{+\infty} \Phi_i u_x \, \varphi_{c_i}(\cdot-x_i) -2
 \int^{x_i}_{-\infty}\Phi_i u_x \, \varphi_{c_i}(\cdot-x_i)  \\
 & =& E_i(u)+E_i (R_X) -2 \int \Phi_i u \, \varphi_{c_i}(\cdot-x_i)
  +2 \int \Phi_i' u \, \varphi_{c_i}(\cdot-x_i)\\
 & & +2 
 \int_{z_i}^{+\infty} \Phi_i u_x \, \varphi_{c_i}(\cdot-x_i) -2
 \int^{z_i}_{-\infty}\Phi_i u_x \, \varphi_{c_i}(\cdot-x_i)  \\
 & = &  E_i(u)+E_i (R_X) -4  c_i u(x_i) \Phi_i(x_i)
  +2 \int \Phi_i' u \, \varphi_{c_i}(\cdot-x_i)\\
 & & -2 
 \int_{x_i}^{+\infty} \Phi_i' u \, \varphi_{c_i}(\cdot-x_i) +2
 \int^{x_i}_{-\infty}\Phi_i' u \, \varphi_{c_i}(\cdot-x_i)  \; .
 \end{eqnarray*}
From \re{de1}-\re{de2}, it is easy to check that $
E_i(R_X)=E(\varphi_{c_i}) +O(e^{-\sqrt{L}/8}) $.
 Since 
 $ C/K = O(L^{-1/2})$
 and, in view of \re{de1},    $ \Phi_i({ x_i})=1+O(e^{-L^{1/2}}) $, it follows that 
 $$
 E_i(u)+E_i(\varphi_{c_i})=E_i(u-R_X)+4c_i M_i +\|u\|_{H^1}^2 O(L^{-1/2}) \; .
 $$
 This yields the result by using that $ E(\varphi_{c_i})=2c_i^2 $. 
%%%%%%%%%%%%%%%%%%%%%%%%%%%%%%%%%
\subsection{End of the proof of Theorem~$\ref{mult-peaks}$}\label{end-proof}
\begin{pro}\label{last}
There exists constants $ C,C'>0 $ independent of $ A $ such that

\begin{equation}\label{oo}
I_{k+1,0}\Bigl( t_0,u(t_0)-R_{X(t_0)}\Bigr) =\sum_{i=k+1}^N
E_i^{t_0}\Bigl( u(t_0)-R_{X(t_0)}\Bigr) \le C (\ep+L^{-{1/4}})
\end{equation}
and 
\begin{equation}\label{oo2}
I_{k+1,0}(t_0)=\sum_{i=k+1}^N
E_i^{t_0}(u(t_0))=\sum_{i=k+1}^N E(\varphi_{c_i})  + O (\ep+L^{-{1/4}})\, .
\end{equation}
with $ |O(x)|\le C' x , \, \forall x\in \R_+^* $.
\end{pro}
{\it Proof. } First it is worth noticing 
 that according to Lemma \ref{eloignement},  $u(t_0)$, $(y_{k+1}(t_0),..,y_{N+1}) $, constructed 
 in \re{defyi}, and $ X(t_0)=(x_1(t_0),..,x_N(t_0)) $, constructed in \re{maxi}, satisfy the hypotheses of Lemma \ref{m-p-lemme}.  Indeed, by construction for $ i\in\{k+1,..,N\} $, $ x_i\in [\tilde{x}_i(t_0)-L/4,\tilde{x}_i(t_0)+L/4, ]\subset 
]y_i(t_0)+L/4, y_{i+1}(t_0)-L/4[  $ and  it is easy to check  that $|u(t_0)|\le O(e^{-\sqrt{L}})+O(\alpha)<3 c_i/4 \le |u(x_i)| $ on $ ]y_i(t_0)-L/4,y_{i+1}(t_0)+L/4[\backslash [\tilde{x}_i(t_0)-L/4,\tilde{x}_i(t_0)+L/4 ]  $ so that 
$$
0<u(t_0,x_i(t_0))=\max_{]y_i(t_0)-L/4,y_{i+1}(t_0)+L/4[} u(t_0) \; .
$$
Therefore, setting $M_i=u(t_0,x_i(t_0)) $, $ \delta_i=c_i-M_i$ and 
 taking the sum over
$i=k+1,..,N $ of \re{eq2m} one gets : $$  \sum_{i=k+1}^{N}\Bigl( M_i E_i^{t_0}(u(t_0)) -
F_i^{t_0}(u(t_0))\Bigr) \ge
-\frac{2}{3}\sum_{i=k+1}^{N}M_i^3+O(L^{-{1/2}})$$ 
Note that by
\re{distH1}  and the continuous embedding of $ H^1(\R) $ into $
L^\infty(\R) $,
 $ M_i=c_i+O(\sqrt{\alpha})+O(e^{-L/8})$, and thus
 \begin{equation}
0<M_{k+1}<\cdot \cdot<M_N \mbox{ and } \delta_i<c_i/2 , \, \forall i\in \{k+1,..,N\} \; .\label{zq}
 \end{equation}
  We set   $
\Delta_0^{t_0} F_i(u)=F_i^{t_0}(u(t_0))-F^0(u(0)) $,
 $ \Delta_0^{t_0} E(u)=E^{t_0}(u(t_0))-E^0(u(0)) $,
  $\Delta_0^{t_0} I_{i,\lambda}(u)=I_{i,\lambda}(t_0,u(t_0))-I_{i,\lambda}(0,u(0)) $.
 Using   the Abel transformation and the monotonicity
estimate~\re{monotonicityestim} (note that $ 0\le 1/M_i\le 2 /c_{k+1} $ for 
 $ i\in \{k+1,..,N\} $), we  get
$$
 \sum_{i=k+1}^{N}M_i \Bigl(\Delta_0^{t_0} E(u) -\frac{1}{M_i}
\Delta_0^{t_0} F(u)\Bigr) =\sum_{i=k+1}^{N}(M_i-M_{i-1})\Delta_0^t
I_{i,1/M_i} \leqslant O(e^{- \sigma_0\sqrt{L}})
$$
and thus 
\begin{equation}\label{e3}
   \sum_{i=k+1}^{N}\Bigl( M_i E_i^{0}(u_0) -
F_i^{0}(u_0))\Bigr) \ge
-\frac{2}{3}\sum_{i=k+1}^{N}M_i^3+O(L^{-{1/2}})\; .
\end{equation}
 By \re{ini}, the
exponential decay of the $ \varphi_{c_i} $'s and the $ \Phi_i
$'s, and the definition
 of $ E_i $ and $ F_i $, it is easy to check that
 \begin{equation}\label{zs}
|E_i^0(u_0)-E(\varphi_{c_i})|+ |F_i^0(u_0)-F (\varphi_{c_i})|\le
O(\varepsilon^2)+O(e^{-\sqrt{L}}), \; \forall i\in\{1,..,N\}\, .
\end{equation}
 Injecting this in  $(\ref{e3})$, taking advantage of \re{zq} and using that 
  $ E(\varphi_{c_i})=2c_i^2 $ and $ F(\varphi_{c_i})=4c_i^3/3 $, we obtain
\begin{eqnarray}
\sum_{i=k+1}^{N}(c_i\delta_i^2-\frac{1}{3}\delta_i^3)
&=&\sum_{i=k+1}^{N}\delta_i^2(c_i-\frac{1}{3}\delta_i)\leqslant
O(\ep^2+L^{-{1/2}}) \nonumber \\
& \Longrightarrow  & \label{e4}
\sum_{i=k+1}^{N} \delta_i^2=
O(\ep^2+L^{-{1/2}}).\label{zqq}
\end{eqnarray}
On the other hand, summing \re{eq3m}  for $ i=k+1, .., N $ one gets  
\begin{equation}\label{lo}
I_{k+1,0}(t_0) - \sum_{i=k+1}^N E(\varphi_{c_i}) = \sum_{i=k+1}^N E_i^{t_0}\Bigl(u(t_0)-R_{X(t_0)}\Bigr)
+4 \sum_{i=k+1}^N c_i \delta_i +O(L^{-1/2}) \, .
\end{equation}
Using \re{zq} and 
the almost monotonicity of $ t\mapsto I_{k+1,0}(t) $, we infer that 
$$
\sum_{i=k+1}^N E_i^{t_0}\Bigl(u(t_0)-R_{X(t_0)}\Bigr)\le  I_{k+1,0}(0)- \sum_{i=k+1}^N E(\varphi_{c_i}) 
+O(\varepsilon +L^{-1/4})
$$
 and \re{zs}-\re{e4}  then yield  \re{oo}.
  Finally,  with \re{oo} at hand, \re{oo2} follows directly from \re{zqq}-\re{lo}. \qed \vspace{2mm}\\
Now, it is crucial to note that  (C-H) is invariant by the change of
unknown $ u(t,x)\mapsto -u(t,-x) $. Therefore setting, for any $ v\in H^1(\R) $,
$$
{\tilde I}_{k,0}(t,v) :=\int_{\R} \Psi(y_{k}(t)-x)[v^2(x)+ v_x^2(x)]\, dx \, ,
$$
with 
$$
y_k(t)=\tilde{x}_k(0)+c_k t/2+L/4\, ,
$$
we infer from Proposition \ref{last}  that
\begin{equation}\label{ooo}
{\tilde I}_{k,0}\Bigl(t_0,u(t_0)-R_{X(t_0)}\Bigr) \le  C (\ep+L^{-{1/4}})\; 
\end{equation}
and 
\begin{equation}\label{ooo2}
{\tilde I}_{k,0}(t_0,u(t_0)) = \sum_{i=1}^k E(\varphi_{c_i}) +O(\ep+L^{-{1/4}})\; .
\end{equation}
Hence,
\begin{eqnarray*}
\tilde{I}_{k,0}(t_0,u(t_0))+I_{k+1,0}(t_0,u(t_0)) & =&
\sum_{i=1}^N E(\varphi_{c_i})+O(\ep+L^{-{1/4}})\\
& = & E(u_0)+O(\ep+L^{-{1/4}})\; .
\end{eqnarray*}
Since $ E(u(t_0))=E(u_0) $ we deduce that
$$ \int_{\R}\Bigl[ 1-\Psi(y_{k}(t_0)-x)-\Psi(x-y_{k+1}(t_0))\Bigr] [u^2(t_0,x)+ u_x^2(t_0,x)]\, dx
 =O(\ep +L^{-{1/4}})\; .
$$
Therefore, since $ | 1-\Psi(y_{k}(t_0)-x)-\Psi(x-y_{k+1}(t_0))|\le O(e^{-\sqrt{L}}) $ for  
 $ x\in \R\backslash ]y_k-L/4, y_{k+1}+L/4[ $ and by the exponentional decay 
  of $ \varphi $, \re{vitesse} and  \re{eloi},
   $$ \int_{y_k-L/4}^{y_{k+1}+L/4} |R_{X(t_0)}|^2+|\partial_x R_{X(t_0)})|^2 \le
  O(e^{-\sqrt{L}/4}) \; , $$ it follows that 
  \begin{equation}\int_{\R}\Bigl[ 1-\Psi(y_{k}(t_0)-\cdot)-\Psi(\cdot-y_{k+1}(t_0))\Bigr] [(u(t_0)-R_{X(t_0})^2
  + (u_x(t_0)-\partial_x R_{X(t_0)})^2]
 =O(\ep +L^{-{1/4}}) \; . \label{oooo}
\end{equation}
Combining \re{oo}, \re{ooo} and \re{oooo} we infer that
$$
E(u(t_0)-R_{X(t_0)}) =O(\ep+L^{-{1/4}})\;
$$
which concludes the proof of \re{e12} since, according to Proposition \ref{last},  $ |O(x)|\le C |x| $ for some constant $ C>0 $ independent of $ A $. This proves \re{ini2} whereas \re{gy} Êfollows from \re{vitesse} and \re{prox}. \vspace{2mm}\\
\noindent
{\bf Acknowledgements}  L.M. would   like to thank  the Oberwolfach
Mathematical center
 where this work was initiated.

\end{document}